\documentclass[11pt]{article}
\usepackage{palatino}
\usepackage{epsfig}
\usepackage{ifthen}
\usepackage{latexsym}
\usepackage{amssymb}
\renewcommand{\baselinestretch}{2.0}

\newtheorem{theorem}{Theorem}[section]
\newtheorem{lemma}[theorem]{Lemma}
\newtheorem{proposition}[theorem]{Proposition}
\newtheorem{corollary}[theorem]{Corollary}
\newtheorem{conjecture}[theorem]{Conjecture}%[section]

\newcommand{\putpicture}[2]{%
\vspace{#2}
\begin{picture}(0,0)
\centerline{\includegraphics{#1}}
\end{picture}
}

\newcommand{\Bx}{\hfill\ensuremath{\Box}}
\newcommand{\cart}{\Box}

\begin{document}
\title{Optimal Pebbling in Products of Graphs}
\author{David S. Herscovici\\%}
Department of Computer Science/IDD\\
CL-AC1 \\
Quinnipiac University \\
275 Mount Carmel Avenue \\
Hamden, CT 06518 \\
\texttt{David.Herscovici@quinnipiac.edu}
\and Benjamin D. Hester \\
Department of Mathematics and Statistics \\
Arizona State University \\
Tempe, AZ, 85287 \\
\texttt{benjamin@mathpost.la.asu.edu}
\and Glenn H. Hurlbert \\
Department of Mathematics and Statistics \\
Arizona State University \\
Tempe, AZ, 85287 \\
\texttt{hurlbert@asu.edu}
}
\renewcommand{\baselinestretch}{1.0}
\maketitle
\renewcommand{\baselinestretch}{2.0}

\begin{abstract}

We prove a generalization of Graham's Conjecture for optimal pebbling
with arbitrary sets of target distributions.  We provide bounds on
optimal pebbling numbers of products of complete graphs and explicitly
find optimal $t$-pebbling numbers for specific such products.  We
obtain bounds on optimal pebbling numbers of powers of the cycle
$C_5$.  Finally, we present explicit distributions which provide
asymptotic bounds on optimal pebbling numbers of hypercubes.

\end{abstract}

{\bf Keywords.} pebbling distribution, pebbling number, fractional pebbling

{\bf 2000 MSC.}  05C99

\newpage

\section{Introduction}
\label{intro}

For a graph $G=(V,E)$, a function $D : V \rightarrow \mathbb{N}$ is
called a \emph{distribution on the vertices of $G$}, or a
\emph{distribution on $G$}.  We usually imagine that $D(v)$ pebbles
are placed on $v$ for each vertex $v \in V$.  Let $|D|$ denote the
\emph{size} of $D$, i.~e.\ $|D| = \displaystyle{\sum_{v \in V} D(v)}$.
For two distributions $D$ and $D'$ on $G$, we say that $D$
\emph{contains} $D'$ if $D'(v)\leq D(v)$ for all $v\in V$.  We allow
pebbling moves on the graph, and define the \emph{pebbling number},
the \emph{optimal pebbling number}, the \emph{$t$-pebbling number},
and the \emph{optimal $t$-pebbling number} of a graph as follows: 

\noindent
\textbf{Definitions}: A \emph{pebbling move} in $G$ takes two pebbles
from a vertex $v \in V$, which contains at least two pebbles, and
places a pebble on a neighbor of $v$.  For two distributions $D_1$ and
$D_2$, we say that $D_2$ is \emph{reachable} from $D_1$ if there is
some sequence of pebbling moves beginning with $D_1$ and resulting in
a distribution which contains $D_2$.  
We say the distribution $D$ is
\emph{solvable}, (respectively, \emph{$t$-solvable}), if every
distribution with one pebble (respectively, $t$ pebbles) on a single
vertex is reachable from $D$. 

The traditional \emph{pebbling number}, and \emph{$t$-pebbling number}
of a graph $G$, denoted $\pi(G)$ and $\pi_t(G)$ respectively, were
defined by Chung~\cite{Hypercubes}.  The \emph{optimal pebbling
number} and \emph{optimal $t$-pebbling number} of $G$, denoted
$\pi^*(G)$ and $\pi_t^*(G)$ respectively, were defined by Pacther,
Snevily, and Voxman~\cite{PSV}.  We give those definitions now. 

\noindent
\textbf{Definitions (Chung~\cite{Hypercubes} and Pachter \emph{et
al.}~\cite{PSV})}: The \emph{$t$-pebbling number} of $G$ is the
smallest number $\pi_t(G)$ such that every distribution $D$ with $|D| \geq
\pi_t(G)$ is $t$-solvable.  The \emph{optimal $t$-pebbling number} of
$G$, denoted $\pi_t^*(G)$, is the smallest number such that some
distribution with $\pi_t(G)$ pebbles is $t$-solvable.  In both cases
we omit the $t$ when $t=1$.  Thus, the \emph{pebbling number} of $G$
is $\pi(G) = \pi_1(G)$ and the \emph{optimal pebbling number} of $G$
is $\pi^*(G) = \pi_1^*(G)$.

The pebbling number was generalized in~\cite{generalizations} to allow
for an arbitrary set of target distributions.  We define this
generalization and extend it to define the optimal pebbling number of
a set of distributions on $G$. 

\noindent
\textbf{Definitions (\cite{generalizations})}: Let $\mathcal{S}$ be a
set of distributions on a graph $G$.  We say a distribution $D$ is
\emph{$\mathcal{S}$-solvable} if every distribution in $\mathcal{S}$
is reachable from $D$.  The \emph{pebbling number of $\mathcal{S}$ in
$G$}, denoted $\pi(G, \mathcal{S})$, is the smallest number such that
every distribution $D$ with $|D| \geq \pi(G, \mathcal{S})$ is
$\mathcal{S}$-solvable.  The \emph{optimal pebbling number of
$\mathcal{S}$ in $G$}, denoted $\pi^*(G, \mathcal{S})$, is the smallest
number such that some distribution $D$ with $|D| = \pi^*(G,
\mathcal{S})$ is $\mathcal{S}$-solvable.

If $\mathcal{S}_t(G)$ consists of all distributions with $t$ pebbles on
a single vertex, we have $\pi(G, \mathcal{S}_1) = \pi(G)$, $\pi^*(G,
\mathcal{S}_1(G)) = \pi^*(G)$, $\pi(G, \mathcal{S}_t(G)) = \pi_t(G)$, and
$\pi^*(G, \mathcal{S}_t(G)) = \pi_t^*(G)$.

\section{Graham's Conjecture and Generalizations in Optimal Pebbling}
\label{Graham}

Graham's Conjecture asserts a bound on the pebbling number of the
Cartesian product of two graphs. 

\noindent
\textbf{Definition}: If $G = (V, E)$ and $G' = (V', E')$ are two
graphs, their Cartesian product is the graph $G\cart G'$ whose vertex
set is the product
\[
V_{G\cart G'} = V \times V' = \{(x, x') : x \in V, x' \in V' \},
\]
and whose edges are given by
\[
E_{G\cart G'} = \{((x, x'), (y, x')) : (x, y) \in E \} \cup \{((x,
x'), (x, y')) : (x', y') \in E' \}.
\]
We also write $G^d$ for the graph $G \cart G \cart \cdots \cart G$
with $d$ copies of $G$ in the product.  Throughout this paper we
follow that convention that $G = (V, E)$ and $G' = (V', E')$.

Chung~\cite{Hypercubes} attributed Conjecture~\ref{Graham's
conjecture} to Graham.
\begin{conjecture}[Graham's Conjecture]
\label{Graham's conjecture}
For any graphs $G$ and $G'$, we have $\pi(G \cart G') \leq \pi(G)
\pi(G')$.
\end{conjecture}
Conjecture~\ref{Graham's conjecture} was generalized
in~\cite{generalizations} to accommodate the more general definitions
of pebbling numbers with arbitrary sets of target distributions.  The
following definition of products of distributions first appeared
in~\cite{all cycles} and the definition of products of sets of
distributions appeared in~\cite{generalizations}.

\noindent
\textbf{Definition (\cite{all cycles, generalizations})}: If $D$ and
$D'$ are distributions on $G$ and $G'$ respectively, then we define
$D \cdot D'$ as the distribution on $G \cart G'$ such that
\[ (D \cdot D') ((x, x')) = D (x) D' (x') \]
for every vertex $(x, x') \in V(G \cart G')$.  Similarly, if
$\mathcal{S}$ and $\mathcal{S}'$ are sets of distributions on $G$
and $G'$ respectively, then $\mathcal{S} \cdot \mathcal{S}'$ is the
set of distributions on $G \cart G'$ given by
\[
\mathcal{S} \cdot \mathcal{S}' = \{ D \cdot D' : D \in
\mathcal{S} \mbox{ and } D' \in \mathcal{S}' \}
\]
Also, for any integer $s$ we define the distribution $s D$ by $(s
D)(x) = s D(x)$ for all $x \in V$.
\begin{conjecture}[\cite{generalizations}]
For all graphs $G$ and $G'$, and all sets of distributions
$\mathcal{S}$ and $\mathcal{S}'$ on $G$ and $G'$ respectively, we
have $\pi(G \cart G', \mathcal{S} \cdot \mathcal{S}') \leq
\pi(G, \mathcal{S}) \pi(G', \mathcal{S}')$.
\label{general Graham}
\end{conjecture}

In this section we prove the analog of Conjecture~\ref{general Graham}
for optimal pebbling.
\begin{theorem}
Let $D$ be an $\mathcal{S}$-solvable distribution on $G$ and let $D'$
be an $\mathcal{S}'$-solvable distribution on $G'$.  Then $D \cdot D'$
is an $(\mathcal{S} \cdot \mathcal{S}')$-solvable distribution on $G
\cart G'$.  In particular, we have $\pi^*(G \cart G', \mathcal{S}
\cdot \mathcal{S}') \leq \pi^*(G, \mathcal{S}) \pi^*(G',
\mathcal{S}')$.
\label{optimal pebbling distributions}
\end{theorem}
To show this, we first establish a few lemmas.
\begin{lemma}
If $D_1$ and $D_2$ are distributions on the graph $G$ such that $D_2$
is reachable from $D_1$, then for any integer $s$, the distribution $s
D_2$ is reachable from $s D_1$.
\label{multiplying by s}
\end{lemma}
\textbf{Proof}: The distribution $s D_1$ may be regarded as $s$
distinct copies of $D_1$.  We can reach $D_2$ from each copy of $D_1$,
so $s D_2$ is reachable from $s D_1$.~\Bx
\begin{lemma}
Let $G$ and $G'$ be graphs.  If $D_1$ and $D_2$ are distributions on
$G$ such that $D_2$ is reachable from $D_1$, then for any distribution
$D'$ on $G'$, $D_2 \cdot D'$ is reachable from $D_1 \cdot D'$.
\label{multiplying by D'}
\end{lemma}
\textbf{Proof}: For each $(x_i, y_j) \in V(G \cart G')$, the number of
pebbles on $(x_i, y_j)$ in the distribution $D_1 \cdot D'$ is given
by $(D_1 \cdot D') ((x_i, y_j)) = D_1(x_i) D'(y_j)$.  Fix $y_j \in
V'$.  We write $G \cart \{ y_j \}$ for the subgraph of $G \cart G'$
induced by the vertices whose second coordinate is $y_j$.  Then $G
\cart \{ y_j \} \cong G$, and if we restrict $D_1 \cdot D'$ to $G
\cart \{ y_j \}$, we obtain the distribution $D'(y_j) D_1$.  Since
$y_j$ is fixed, $D'(y_j)$ is a constant, so by Lemma~\ref{multiplying
by s}, the distribution $D'(y_j) D_2$ is reachable in $G \cart \{
y_j \}$.  Repeating this for each $y_j \in V'$, we end up with a
distribution in which each $(x_i, y_j)$ has at least $D_2(x_i)
D'(y_j) = (D_2 \cdot D') ((x_i, y_j))$ pebbles, so $D_2 \cdot D'$
is reachable from $D_1 \cdot D'$.~\Bx 

We are now ready to prove Theorem~\ref{optimal pebbling
distributions}. 

\noindent
\textbf{Proof of Theorem~\ref{optimal pebbling distributions}}: Let
$D$ and $D'$ be $\mathcal{S}$- and $\mathcal{S}'$-solvable
distributions on $G$ and $G'$ respectively.  To show that $D \cdot D'$
is $(\mathcal{S} \cdot \mathcal{S}')$-solvable on $G\cart G'$, let
$\Delta$ be a distribution in $\mathcal{S} \cdot \mathcal{S}'$.  Then
we can write $\Delta = D_i \cdot D_j'$ for some $D_i \in \mathcal{S}$
and $D_j' \in \mathcal{S}'$.  Also, $D_i$ is reachable from $D$ and
$D_j'$ is reachable from $D'$.  Thus, by Lemma~\ref{multiplying by
D'}, $D_i \cdot D_j'$ is reachable from $D \cdot D_j'$, which is
reachable from $D \cdot D'$.

If we choose $D$ and $D'$ so that $|D| = \pi^*(G)$ and $|D'| =
\pi^*(G')$, we have
\[
|D \cdot D'| = \sum_{x_i \in V} \sum_{y_j \in V'} D (x_i) D' (y_j) =
\sum_{x_i \in V} D(x_i) \sum_{y_j \in V'} D' (y_j) = |D| |D'|.
\]
Thus, $D \cdot D'$ is an $(\mathcal{S} \cdot \mathcal{S}')$-solvable
distribution on $G \cart G'$ with $|D| |D'| = \pi^*(G, \mathcal{S})
\pi^*(G', \mathcal{S})$; therefore $\pi^*(G \cart G', \mathcal{S}
\cdot \mathcal{S}') \leq \pi^*(G, \mathcal{S}) \pi^*(G',
\mathcal{S}')$, as desired.~\Bx 

Corollaries~\ref{optimal st-pebbling graphs} and~\ref{optimal pebbling
graphs} follow immediately from Theorem~\ref{optimal pebbling
distributions}.  Fu and Shiue~\cite{Fu Shiue} announced
Corollary~\ref{optimal pebbling graphs}, the optimal pebbling analog to
Graham's Conjecture.  Shiue proved it in~\cite{Shiue}.
\begin{corollary}
For all graphs $G$ and $G'$ and all positive integers $s$ and $t$, we
have $\pi_{st}^*(G \cart G') \leq \pi_s^*(G) \pi_t^*(G')$.
\label{optimal st-pebbling graphs}
\end{corollary}
\begin{corollary}[Fu and Shiue~\cite{Fu Shiue, Shiue}]
For all graphs $G$ and $G'$, we have $\pi^*(G \cart G') \leq \pi^*(G)
\pi^*(G')$.
\label{optimal pebbling graphs}
\end{corollary}

\section{Products of Complete Graphs}
\label{kn products section}

Our work in Section~\ref{Graham} puts an upper bound on $\pi_t^*(G
\cart G')$.  In this section, we improve those bounds when $G$ and
$G'$ are complete graphs.  Our main result is Theorem~\ref{bounds for
pi^*(GxKn)}.
\begin{theorem}
For any graph $G$ and any positive integer $t$, we have $\left\lceil
\left( \frac{n}{n+1} \right) \pi^*_{2t}(G) \right\rceil \leq \pi_t^*(G
\cart K_n) \leq \pi^*_{2t} (G)$.
\label{bounds for pi^*(GxKn)}
\end{theorem}
Lemma~\ref{at most one odd vertex} helps us find the optimal
$t$-pebbling number of a complete graph (Theorem~\ref{pi t(Kn)}).

\noindent
\textbf{Definition}: Given any distribution of pebbles on the vertices
of the graph $G$, we say the vertex $v$ is \emph{odd} or \emph{even},
depending on whether it has an odd or an even number of pebbles on it.
\begin{lemma}
Let $t$ and $n$ be positive integers, and suppose we have a
$t$-solvable distribution with $\pi_t^*(K_n)$ pebbles on the vertices
of $K_n$.  Then:
\begin{enumerate}
\item If the vertex $v_i$ is odd, then every other vertex has at least
as many pebbles as $v_i$.
\item There are at most two odd vertices.
\item If there are two odd vertices in $K_n$, then moving a pebble
from one of these vertices to the other creates another $t$-solvable
distribution.
\end{enumerate}
In particular, some $t$-solvable distribution of $\pi_t^*(K_n)$ pebbles
on $K_n$ has at most one odd vertex.
\label{at most one odd vertex}
\end{lemma}
\textbf{Proof}: Removing a pebble from an odd vertex $v_i$ does not
affect the number of pebbles that may be moved to any other vertex;
thus, every other vertex may still receive $t$ pebbles.  Since there
would now be fewer than $\pi_t^*(K_n)$ pebbles, $v_i$ could no longer
receive $t$ pebbles.  If another vertex $v_j$ started with fewer
pebbles than $v_i$, we could use the pebbles now on $v_i$ and $v_j$ to
put at least as many pebbles on $v_i$ as on $v_j$, and any pebbles
that could be moved to $v_j$ from other vertices could also be moved
to $v_i$.  Thus, we could put at least as many pebbles on $v_i$ as on
$v_j$, contradicting our assertion that $t$ pebbles can be moved to
$v_j$, but not to $v_i$.  Therefore, every other vertex has at least
as many pebbles as the odd vertex $v_i$.

If there are two or more odd vertices in $K_n$, we remove a pebble
from each of these vertices and add two pebbles to any vertex, say
$v_1$.  Now every vertex can receive at least as many pebbles as it
could from the original distribution: if the target originally was
odd, the first move would be from $v_1$ to the target.  We therefore
have a $t$-solvable distribution in which every vertex is even.
Furthermore, if we originally had three or more odd vertices, this
distribution would have fewer pebbles, contradicting the hypothesis
that the original distribution had $\pi_t^*(K_n)$ pebbles.~\Bx
\begin{theorem}
For any positive integers $n$ and $t$, let $q = t \mbox{
\textnormal{div} } (n+1)$ and let $r = t \bmod{(n+1)}$.  Thus, $t =
(n+1)q + r$.  Then $\pi_t^*(K_n)$ is given by
\[
\pi_t^*(K_n) = \left\{ \begin{array}{ll}
2t - 2q = 2nq + 2r & \mbox{if } r < n \\
2t - 2q - 1 = 2nq + 2n - 1 & \mbox{if } r = n \\
\end{array}\right.
\]
In particular, $\pi_t^*(K_n) = 2t$ if and only if $t < n$.
\label{pi t(Kn)}
\end{theorem}
\textbf{Proof}: First note that if we put $2q+2r$ pebbles on one
vertex and we put $2q$ pebbles on every other vertex, then we can move
an additional $(n-1)q + r$ pebbles onto any vertex that starts with
$2q$ pebbles, and we can move $(n-1)q$ additional pebbles onto the
vertex that starts with $2q+2r$ pebbles.  In either case, we can move
at least $t = (n+1)q + r$ pebbles to any target, including the pebbles
that start there.  Thus, $\pi_t^*(K_n) \leq 2nq + 2r$.

We now consider whether a $t$-solvable distribution in $K_n$ could
have fewer than $2nq + 2r$ pebbles.  Let $v_i$ be the vertex with the
fewest pebbles, and suppose it has $p_i$ pebbles.  Adding $t-p_i$
pebbles to $v_i$ costs at least $2(t-p_i)$ pebbles.  Therefore,
including the pebbles that started on $v_i$, the original distribution
has at least $2t - p_i$ pebbles.  If this is less than $2nq + 2r = 2t
- 2q$, then $p_i > 2q$.

If $p_i \geq 2q+2$, every vertex has at least $2q+2$ pebbles, and so
the distribution uses $(2q+2)n = 2nq + 2n \geq 2nq + 2r$ pebbles.
Therefore, we assume $p_i = 2q+1$.  Now by Lemma~\ref{at most one odd
vertex}, we may assume every other vertex has at least $2q+2$ pebbles,
so we have already accounted for $(2q+2) n - 1 = 2nq + 2n - 1$
pebbles.  The only way this can be smaller than $2nq + 2r$ is if $r =
n$.  Now we simply observe that if $r = n$ and $2q+1$ pebbles are on
$v_i$ and $2q+2$ pebbles are on every other vertex, then a total of
$(2q+1) + (n-1) (q+1) = (n+1)q + n = t$ pebbles can be moved to $v_i$,
and similarly, $(2q+2) + (n-2) (q+1) + q = (n+1)q + n = t$ pebbles can
be move to any other vertex.  Finally, $\pi_t^*(K_n) = 2t$ if and only if
$q = 0$ and $r < n$, i.~e.\ if and only if $t < n$.~\Bx

The optimal $t$-pebbling number is not generally monotone, in the
following sense.  If it is large for a particular graph, it can be
reduced significantly by the addition of a single vertex adjacent to
all others.  However, for complete graphs the parameter is
nondecreasing.
\begin{proposition}
For every graph $G$ and every positive integer $n$, we have $\pi_t^*(K_n
\cart G) \leq \pi_t^*(K_{n+1} \cart G)$.
\label{K(n+1) x G}
\end{proposition}
\textbf{Proof}: Given any distribution $D: G \cart K_{n+1} \rightarrow
\mathbb{N}$, let $g(D): G \cart K_n \rightarrow \mathbb{N}$ be the
distribution on $G \cart K_n$ defined by
\[ (g(D))(v, w_i) = \left \{ \begin{array}{ll}
D(v, w_i) & \mbox{ if $i < n$} \\
D(v, w_n) + D(v, w_{n+1}) & \mbox{if $i = n$}.
\end{array} \right. \]
Then any pebbling move from $D$ to $D'$ in $G \cart K_{n+1}$ can be
shadowed by moves from $g(D)$ to a distribution that contains $g(D')$
in $G \cart K_n$: moves from $(v, w_n)$ to $(v, w_{n+1})$ or
\emph{vice versa} may be ignored, other moves from $D$ to $D'$ either
from, to, or within $G \cart \{ w_{n+1} \}$ can be made from $g(D)$ to
$g(D')$ using $G \cart \{ v_n \}$ instead, and moves from $D$ to $D'$
that do not use $G \cart \{ v_{n+1} \}$ can be made unchanged from
$g(D)$ to $g(D')$.  Therefore, if $D$ is a $t$-solvable distribution
on $G \cart K_{n+1}$ then $g(D)$ is a $t$-solvable distribution on
$G \cart K_n$.  Since $|g(D)| = |D|$, we have $\pi_t^*(G \cart K_n)
\leq \pi_t^*(G \cart K_{n+1})$.~\Bx

Corollary~\ref{pi(G x K_m) <= pi(G x K_n)} follows from
Proposition~\ref{K(n+1) x G} by induction on $n$, starting with $n=m$
as a basis.
\begin{corollary}
For every graph $G$ and all positive integers $m$ and $n$ with $m \leq
n$, we have $\pi_t^*(G \cart K_m) \leq \pi_t^*(G \cart K_n)$.\Bx
\label{pi(G x K_m) <= pi(G x K_n)}
\end{corollary}
\textbf{Definitions}: Given a distribution $D:V(G \cart G')
\rightarrow \mathbb{N}$ on $G \cart G'$ and a subset $S \subseteq
V'$, we define the distribution $f_S(D) : V \rightarrow
\mathbb{N}$ on $G$ by
\[
(f_S(D))(u) = \sum_{v \in V' \setminus S} D(u, v) + 2 \sum_{v
\in S} D(u, v).
\]
for every $u \in V$.  In other words, we count every pebble on a
vertex whose coordinate in $G'$ is in $S$ twice and every other pebble
once.  If the vertices of $G'$ are $\{ v_1, v_2, \ldots, v_n \}$, we
define $f_i(D)$ by
\[
(f_i(D))(u) = (f_{\{v_i\}}(D))(u) = \sum_{j \neq i} D(u, v_j) + 2
D(u, v_i).
\]

Lemmas~\ref{matching sequence from pi(D)} and~\ref{useless pebble} are
key to proving Theorem~\ref{G x Kn}, which is the upper bound in
Theorem~\ref{bounds for pi^*(GxKn)}.
\begin{lemma}
Let $S$ be any nonempty subset $S \subseteq V'$, and suppose there
is a sequence of pebbling moves in $G \cart G'$ from $D_0$ to $D_k$.
Then there is a sequence of pebbling moves in $G$ from $f_S(D_0)$ to
a distribution that contains $f_S(D_k)$.  In particular, if
$|f_S(D_0)| < \pi^*_{2t}(G)$, then $D_0$ cannot be $t$-solvable in $G
\cart G'$.
\label{matching sequence from pi(D)}
\end{lemma}
\textbf{Proof}: Let $D_0, D_1, \ldots, D_k$ be the sequence of
distributions in $G \cart G'$ after each pebbling move.  We show by
induction that we can shadow each pebbling move in $G \cart G'$
with moves in $G$.  Toward that end, suppose that there is a sequence
of pebbling moves in $G$ from $f_S(D_0)$ to a distribution that
contains $f_S(D_i)$.  The basis $i = 0$ is trivial.

Suppose going from $D_i$ to $D_{i+1}$ requires a move from $(u, v_1)$
to $(u, v_2)$.  Then the pebbles involved in the move add either four
or two pebbles to $u$ in $f_S (D_i)$, depending on whether $v_1 \in
S$, and they add either two pebbles or one pebble to $u$ in $f_S
(D_{i+1})$, depending on whether $v_2 \in S$ or not.  In either case,
$f_S (D_i)$ contains $f_S (D_{i+1})$, and we can simply ignore the
extra pebbles.

Otherwise, going from $D_i$ to $D_{i+1}$ requires a move from $(u_1,
v)$ to $(u_2, v)$.  If $v \in S$ the pebbles involved in this move add
four pebbles to $u_1$ and two pebbles to $u_2$ in $f_S(D_i)$ and
$f_S(D_{i+1})$, and if $v \notin S$, they add two pebbles to $u_1$
and one pebble to $u_2$ in $f_S(D_i)$ and $f_S(D_{i+1})$,
respectively.  The latter case simply requires a pebbling move from
$u_1$ to $u_2$ in $G$ to get from $f_S(D_i)$ to $f_S(D_{i+1})$;
the former case requires two such moves.

In any of these cases, we can go from $f_S(D_0)$ to a distribution
that contains $f_S(D_i)$ to one that contains $f_S(D_{i+1})$.
Continuing this process, we reach a distribution that contains
$f_S(D_k)$.

Now if $|f_S(D_0)| < \pi^*_{2t}(G)$, there is some vertex $x \in V$
such that $2t$ pebbles cannot be moved onto $x$ by any sequence of
pebbling moves starting from $f_S(D_0)$.  Therefore, we cannot reach
any distribution $D_k$ in $G \cart G'$ for which $(f_S(D_k))(x) \geq
2t$.  In particular, for any $s \in S$, we cannot move $t$ pebbles
onto the vertex $(x, s)$.~\Bx

Lemma~\ref{useless pebble} tells us that if some copy of $G$ in $G
\cart K_n$ starts with a single pebble, then that pebble does not
help us reach vertices in any other copy of $G$.
\begin{lemma}
Let $D:V(G \cart K_n) \rightarrow \mathbb{N}$ be a distribution of
pebbles on $G \cart K_n$, and suppose there is at most one pebble on
some $G \cart \{v_i\}$.  Let $D'$ be the distribution on $G \cart K_n$
obtained by removing that pebble, or let $D' = D$ if there is no such
pebble.  Let $S = V(K_{n}) \setminus \{ v_i \}$, and let $D''$ be any
configuration of pebbles on $G \cart S$ that we can reach from $D$.
Then we can reach a configuration that contains $D''$ starting from
$D'$.
\label{useless pebble}
\end{lemma}
\textbf{Proof}: If there are no pebbles on $G \cart \{ v_i \}$ and $D
= D'$, there is nothing to prove, so we assume there is a pebble on $G
\cart \{ v_i \}$ in $D$.  Paint this pebble gold, and assume it
survives every pebbling move in the sequence from $D$ to $D''$ in
which it participates.

If the gold pebble never leaves $G \cart \{ v_i \}$, we can make the
same moves in $D'$ as in $D$ and ignore the moves involving the gold
pebble.  Otherwise, let $v_j$ be the vertex in $K_n$ involved in the
first move of the gold pebble from $(x, v_i)$ to $(x, v_j)$.  We
examine the moves by the gold pebble before it leaves $G \cart
\{ v_i \}$.  Note that every such move consumes a nongold pebble that
was moved onto $G \cart \{ v_i \}$ from a different copy of $G$.  Our
approach is to move those pebbles to $G \cart \{ v_j \}$ instead.

Thus, from $D'$, we ignore all moves from $D$ involving the gold
pebble before it first leaves $G \cart \{ v_i \}$.  We replace all
other moves to, from, or within $G \cart \{ v_i \}$ with moves to,
from, or within $G \cart \{ v_j \}$, ignoring moves between $G \cart
\{ v_i \}$ and $G \cart \{ v_j \}$.  Now the pebble that would have
been removed from $(x, v_i)$ when the gold pebble moved to $(x, v_j)$
reaches $(x, v_j)$ in place of the gold pebble.  This pebble can
replace of the gold pebble on all subsequent moves.  The result of
these changes is that all pebbles that ended up on $G \cart S$
starting from $D$ end up on the same vertices starting from $D'$,
except that the gold pebble is replaced by a different pebble.~\Bx 

\noindent
\textbf{Notation}: Suppose we have a distribution of pebbles on $G
\cart K_n$.  For each $i$ with $1 \leq i \leq n$, we let $p_i$ be the
number of pebbles on $G \cart \{ v_i \}$, and we assume without loss
of generality that $p_1 \leq p_2 \leq \cdots \leq p_n$.

Theorem~\ref{G x Kn} gives the upper bound from Theorem~\ref{bounds
for pi^*(GxKn)}.
\begin{theorem}
For any graph $G$ and any positive integer $n$, we have $\pi_t^*(G \cart
K_n) \leq \pi^*_{2t} (G)$.  Furthermore, equality holds when $2n \geq
\pi^*_{2t}(G) + 1$.
\label{G x Kn}
\end{theorem}
\textbf{Proof}: We first note that if $D$ is a $2t$-solvable
distribution on $G$, then placing $D(x)$ pebbles on the vertex $(x,
v_1)$ for every $x \in V$ creates a distribution from which $t$
pebbles can be moved to the vertex $(x_i, v_j)$ since we can first
move $2t$ pebbles to $(x_i, v_1)$.  Therefore, $\pi_t^*(G \cart K_n)
\leq \pi^*_{2t} (G)$.

Now suppose $2n \geq \pi^*_{2t}(G) + 1$, and let $D$ be a distribution
on $G \cart K_n$ with $\pi^*_{2t}(G) - 1$ pebbles or fewer.  Then
either $p_1 = 0$ or $p_1 = p_2 = 1$; otherwise, we would have $1 \leq
p_1$ and $2 \leq p_2 \leq p_3 \leq \cdots \leq p_n$.  But then $|D|
\geq 2n-1 \geq \pi^*_{2t}(G)$, contrary to our assumption that $D$ has
at most $\pi^*_{2t}(G) - 1$ pebbles.

If $p_1 = 0$, then $f_1 (D)$ has at most $\pi^*_{2t}(G) - 1$ pebbles,
so $2t$ pebbles cannot be moved onto some $x \in V$ starting from $f_1
(D)$.  Therefore, by Lemma~\ref{matching sequence from pi(D)}, we
cannot move $t$ pebbles onto $(x, v_1)$ starting from $D$.  On the
other hand, if $p_1 = p_2 = 1$, let $D'$ be the distribution on $G
\cart K_n$ with the lone pebble on $G \cart \{ v_2 \}$ removed.
Then $|f_1(D')| \leq \pi^*_{2t}(G) - 1$, since the pebble on $G \cart
\{ v_1 \}$ that is counted twice is offset by the pebble that is
removed from $G \cart \{ v_2 \}$.  As before, Theorem~\ref{matching
sequence from pi(D)}, shows that $t$ pebbles cannot be moved to some
$(x, v_1)$ in $V(G \cart K_n)$ starting from $D'$.  But now applying
Lemma~\ref{useless pebble} with $i = 2$ shows that $t$ pebbles
cannot be moved to $(x, v_1)$ from $D$ in this case either.
Therefore, $\pi_t^*(G \cart K_n) = \pi^*_{2t} (G)$.~\Bx

Applying Theorem~\ref{G x Kn} inductively gives Corollary~\ref{G x
K_n1 x K_n2 x ... x K_nd}.
\begin{corollary}
For any graph $G$, any positive integer $t$, and any sequence of
integers $n_1, n_2, \ldots n_d$, we have
\[
\pi_t^*(G \cart K_{n_1} \cart K_{n_2} \cart \cdots \cart K_{n_d})
\leq \pi^*_{2^d t}(G).
\]
Furthermore, equality holds if $2n_i \geq \pi^*_{2^d t}(G)+1$ for each
$n_i$.
\label{G x K_n1 x K_n2 x ... x K_nd}
\end{corollary}
\textbf{Proof}: We fix $d$ and $t$, and prove by induction on $k$ that
\begin{equation}
\pi^*_{2^{d-k} t}(G \cart K_{n_1} \cart K_{n_2} \cart \cdots \cart
K_{n_k}) \leq \pi^*_{2^d t}(G),
\label{G x K_n1 x ... x K_nd equation}
\end{equation}
and that equality holds when each $n_i$ satisfies $2n_i \geq
\pi^*_{2^d t}(G)+1$.  The basis $k = 0$ is trivial, so we assume
that~(\ref{G x K_n1 x ... x K_nd equation}) holds for some $k$ with $0
\leq k < d$.  Applying Theorem~\ref{G x Kn} and then applying~(\ref{G
x K_n1 x ... x K_nd equation}) gives
\[
\pi^*_{2^{d-k-1} t} (G \cart K_{n_1} \cart K_{n_2} \cart \cdots
\cart K_{n_k} \cart K_{n_{k+1}}) \leq \pi^*_{2^{d-k} t}(G \cart
K_{n_1} \cart K_{n_2} \cart \cdots \cart K_{n_k}) \leq \pi^*_{2^d
t}(G).
\]
as desired.  Furthermore, equality continues to hold if $n_{k+1}$
satisfies $2 n_{k+1} \geq \pi^*_{2^d t}(G)+1$.~\Bx
\begin{corollary}
For all positive integers $t$, and any product of $d$ complete graphs,
we have
\[ \pi_t^*(K_{n_1} \cart K_{n_2} \cart \cdots \cart K_{n_d}) = 2^d t
\]
if and only if each $n_i \geq 2^{d-1} t + 1$.
\label{kn products}
\end{corollary}
\textbf{Proof}: Applying Corollary~\ref{G x K_n1 x K_n2 x ... x K_nd}
with $G$ equal to the trivial graph gives
$\pi_t^* (K_{n_1} \cart K_{n_2} \cart \cdots \cart K_{n_d}) \leq
\pi^*_{2^d t} (G) = 2^d t$.
Furthermore, equality holds when each $n_i$ satisfies $2 n_i \geq 2^d
t + 1$, or equivalently, $n_i \geq 2^{d-1} t + 1$.  On the other hand,
if $n_i \leq 2^{d-1} t$ for some $i$, we assume without loss of
generality that $n_1 \leq 2^{d-1} t$.  Now applying Corollary~\ref{G x
K_n1 x K_n2 x ... x K_nd} with $G = K_{n_1}$ gives
$\pi_t^* (K_{n_1} \cart K_{n_2} \cart \cdots \cart K_{n_d}) \leq
\pi^*_{2^{d-1}t} (K_{n_1})$,
and by Proposition~\ref{pi t(Kn)}, $\pi^*_{2^{d-1}t} (K_{n_1}) \leq 2^d
t - 1$ when $n_1 \leq 2^{d-1} t$.~\Bx

We can now prove Theorem~\ref{bounds for pi^*(GxKn)}.

\noindent
\textbf{Proof of Theorem \ref{bounds for pi^*(GxKn)}}: 
The upper bound is given by Theorem~\ref{G x Kn}.  To
establish the lower bound, suppose we have a $t$-solvable distribution
$D$ of $P = \pi_t^* (G \cart K_n)$ pebbles on $G \cart K_n$.  Since
$p_1 \leq p_2 \leq \cdots \leq p_n$, we have $p_1 \leq \frac{P}{n}$.
Now by Lemma~\ref{matching sequence from pi(D)}, $\pi^*_{2t}(G) \leq
|f_1(D)| = P + p_1 \leq P + \frac{P}{n} = \left(\frac{n+1}{n} \right)
P$.  Since $P$ must be an integer, we have $P \geq \left\lceil \left(
\frac{n}{n+1} \right) \pi^*_{2t}(G) \right\rceil$.~\Bx

For the smallest of products, we are able to get exact results for all $t$.
These exhibit a nice pattern that we will say more about subsequently.
First we present an obvious proposition.
\begin{proposition}
For any graph $G$ and any positive integers $s$ and $t$, we have
$\pi^*_{s+t}(G) \leq \pi^*_s(G) + \pi_t^*(G)$.  Similarly, for regular
pebbling, we have $\pi_{s+t}(G) \leq \pi_s(G) + \pi_t(G)$.
\label{pi^*_s+t}
\end{proposition}
\textbf{Proof}: We can place $\pi^*_s(G)$ red pebbles and $\pi_t^*(G)$
blue pebbles on $G$ in such a way that $s$ red pebbles and $t$ blue
pebbles can be moved to any target vertex.

For regular pebbling, we note that from any placement of $\pi_s(G) +
\pi_t(G)$ pebbles, if we arbitrarily paint $\pi_s(G)$ pebbles red and
$\pi_t(G)$ pebbles blue, then $s$ red pebbles and $t$ blue pebbles
can be moved to any target vertex.~\Bx
\begin{proposition}
To find the optimal $t$-pebbling number of $K_2 \cart K_2$, let $q =
t \mbox{ div } 9$ and $r = t \bmod{9}$.  Then
\[
\pi_t^*(K_2 \cart K_2) = \left\{ \begin{array}{ll}
3 & \mbox{if } t=1 \\
16q + 2r & \mbox{if } r \in \{ 0, 1, 2, 3, 4, 5 \} \mbox{ and } t \neq 1 \\
16q + 2r - 1 & \mbox{if } r \in \{ 6, 7, 8 \}. \\
\end{array}\right.
\]
In each case except $t = 1$, the lower bound from Theorem~\ref{bounds
for pi^*(GxKn)} is tight.
\label{pi t(K2xK2)}
\end{proposition}
\textbf{Proof}: If $t=1$, we note that two pebbles are not enough to
reach every vertex: if we put them on different vertices, the
unoccupied vertices cannot be reached, and if we put them on the same
vertex, the antipodal vertex is unreachable.  On the other hand, three
pebbles are sufficient, since we can put two pebbles on $(v_0, v_0)$
and one on $(v_1, v_1)$.

For $2 \leq t \leq 10$, we consider Table~\ref{K2 x K2 table}: 
\begin{table}[ht]
\[
\begin{array}{|*{10}{c|}}
\hline
\mathbf{t} & 2 & 3 & 4 & 5 & 6 & 7 & 8 & 9 & 10 \\
\hline
\mathbf{\left\lceil \frac{2}{3} \pi^*_{2t}(K_2) \right\rceil} &
4 & 6 & 8 & 10 & 11 & 13 & 15 & 16 & 18 \\
\hline
\mbox{\textbf{Optimal Distribution}}
& 2, 0 & 2, 1 & 2, 2 & 3, 2 & 4, 2 & 4, 3 & 4, 4 & 4, 4 & 5, 4 \\
p_{00}, p_{01}/p_{10}, p_{11}
& 0, 2 & 1, 2 & 2, 2 & 2, 3 & 2, 3 & 3, 3 & 3, 4 & 4, 4 & 4, 5 \\
\hline
\end{array}
\]
\caption{Computing $\pi_t^*(K_2 \cart K_2)$ for $2 \leq t \leq 10$}
\label{K2 x K2 table}
\end{table}
The second row of this table gives the lower bound for $\pi_t^* (K_2
\cart K_2)$ from Theorem~\ref{bounds for pi^*(GxKn)}, and the last
row gives a solvable distribution with the given number of pebbles.
Therefore, the bound is tight.

Finally, for $t \geq 11$, we assume by induction on $t$ that the lower
bound is tight for $t' = t-9$, and we show that $\pi_t^*(K_2 \cart
K_2) = \pi^*_{t'}(K_2 \cart K_2) + 16$.  Comparing the computation
of the lower bound for $\pi_t^*(K_2 \cart K_2)$ to that of
$\pi^*_{t'}(K_2 \cart K_2)$, we have $2t = 2t' + 18$, so $2t
\mbox{ div } 3 = 2t' \mbox{ div } 3 + 6$, and $2t \bmod{3} = 2t'
\bmod{3}$.  Thus, $\pi^*_{2t}(K_2) = \pi^*_{2t'}(K_2) + 24$, and
the lower bound from Theorem~\ref{bounds for pi^*(GxKn)} gives
$\pi_t^* (K_2 \cart K_2) \geq \pi^*_{t'} (K_2 \cart K_2) +
16$.  On the other hand, Proposition~\ref{pi^*_s+t} tells us that
$\pi_t^*(K_2 \cart K_2) \leq \pi^*_{t'} (K_2 \cart K_2) +
\pi^*_9 (K_2 \cart K_2) = \pi^*_{t'} (K_2 \cart K_2) + 16$.
Therefore, $\pi_t^* (K_2 \cart K_2) = \pi^*_{t'} (K_2 \cart
K_2) + 16$, as required.~\Bx

We can compute $\pi_t^* (K_2 \cart K_3)$ similarly.
\begin{proposition}
The optimal $t$-pebbling number of $K_2 \cart K_3$ is
\[
\pi_t^*(K_2 \cart K_3) = \max \left( \left\lceil \frac{2}{3}
\pi^*_{2t}(K_3) \right\rceil, \left\lceil \frac{3}{4} \pi^*_{2t}(K_2)
\right\rceil \right).
\]
In particular, if $q = t \mbox{ div } 6$ and $r = t \bmod{6}$, then
\[
\pi_t^*(K_2 \cart K_3) = \left\{ \begin{array}{ll}
12q & \mbox{if } r = 0 \\
12q + 2r + 1 & \mbox{otherwise}. \\
\end{array}\right.
\]
\label{pi t(K2xK3)}
\end{proposition}
\textbf{Proof}: For $1 \leq t \leq 6$, we use Table~\ref{K2 x K3 table}.
\begin{table}[ht]
\[
\begin{array}{|*{7}{c|}}
\hline
\mathbf{t} & 1 & 2 & 3 & 4 & 5 & 6 \\
\hline
 \mathbf{\left\lceil \frac{3}{4} \pi^*_{2t}(K_2) \right\rceil} &
3 & 5 & 6 & 9 & 11 & 12 \\
\hline
 \mathbf{\left\lceil \frac{2}{3} \pi^*_{2t}(K_3) \right\rceil} &
3 & 4 & 7 & 8 & 11 & 12 \\
\hline
\mbox{\textbf{Optimal Distribution}}
& 2, 0, 0 &  2, 0, 1 &  2, 0, 2 &  2, 2, 2 &  2, 2, 2 &  2, 2, 2, \\
p_{00}, p_{01}, p_{02}/p_{10}, p_{11}, p_{12}
& 1, 0, 0 &  0, 2, 0 &  1, 2, 0 &  1, 1, 1 &  2, 2, 1 &  2, 2, 2 \\
\hline
\end{array}
\]
\caption{Computing $\pi_t^*(K_2 \cart K_3)$ for $1 \leq t \leq 6$}
\label{K2 x K3 table}
\end{table}
%
% \begin{table}[ht]
% \[
% \begin{array}{c|c|c|c}
% \mathbf{t} &
%  \mathbf{\left\lceil \frac{3}{4} \pi^*_{2t}(K_2) \right\rceil} &
%  \mathbf{\left\lceil \frac{2}{3} \pi^*_{2t}(K_3) \right\rceil} &
% \mbox{\textbf{OSD}} \\
% &&& (p_{00}, p_{01}, p_{02}; p_{10}, p_{11}, p_{12}) \\
% \hline
% 1 & 3  & 3  & (2, 0, 0; 1, 0, 0) \\
% 2 & 5  & 4  & (2, 0, 1; 0, 2, 0) \\
% 3 & 6  & 7  & (2, 0, 2; 1, 2, 0) \\
% 4 & 9  & 8  & (2, 2, 2; 1, 1, 1) \\
% 5 & 11 & 11 & (2, 2, 2; 2, 2, 1) \\
% 6 & 12 & 12 & (2, 2, 2; 2, 2, 2) \\
% \end{array}
% \]
% \caption{Computing $\pi_t^*(K_2 \cart K_3)$ for $1 \leq t \leq 6$}
% \label{orig K2 x K3 table}
% \end{table}
% 
For larger $t$, we note from Proposition~\ref{pi^*_s+t} that
$\pi^*_{t'+6} (K_2 \cart K_3) \leq \pi^*_{t'} (K_2 \cart K_3) +
\pi^*_6 (K_2 \cart K_3) = \pi^*_{t'} (K_2 \cart K_3) + 12$, which
agrees with the asserted lower bound.~\Bx

Corollary~\ref{kn products} shows that for small values of $t$, the
upper bound in Theorem~\ref{bounds for pi^*(GxKn)} is tight for
products of complete graphs.  It was obtained by applying
Theorem~\ref{G x Kn} inductively, with $G$ being the trivial graph.
If we apply the lower bound in Theorem~\ref{bounds for pi^*(GxKn)}
inductively with $G$ being the trivial graph, we get a lower bound on
the optimal $t$-pebbling number of a product of complete graphs.
Theorem~\ref{asymptotic bound for complete products} shows that this
lower bound is asymptotically tight as $t$ gets large.  We begin with
Lemma~\ref{2^d pebbles on each vertex}
\begin{lemma}
\label{2^d pebbles on each vertex}
Let $n_1, n_2, \ldots, n_d$ be a sequence of nonnegative integers, and
let $T_j = \displaystyle{\prod_{i=1}^j (n_i + 1)}$.  Then for any
integer $k$, putting $2^d k$ pebbles on each vertex of $G = K_{n_1}
\cart K_{n_2} \cart \cdots \cart K_{n_d}$ creates a $kT_d$-solvable
configuration.  Thus, $\pi_{kT_d}^* (G) = 2^d k
\displaystyle{\prod_{i=1}^d n_i}$.
\end{lemma}
\textbf{Proof}: If $d=0$ the products are all empty, so $T_0=1$ and
$G$ is the trivial graph.  Clearly, putting $k$ pebbles on the lone
vertex gives an optimal $k$-solvable configuration, as required.  For
larger $d$, we first show the specified configuration is
$kT_d$-solvable.  Toward that end, Let $(x_1, x_2, \ldots, x_d)$ be
the target vertex in $K_{n_1} \cart K_{n_2} \cart \cdots \cart
K_{n_d}$.  If we have $2^d k$ pebbles on each vertex, then for each $v
\in V(K_{n_d})$, we have $2^{d-1} (2k)$ pebbles on each vertex of
$K_{n_1} \cart K_{n_2} \cart \cdots \cart K_{n_{d-1}} \cart \{ v \}
\cong K_{n_1} \cart K_{n_2} \cart \cdots \cart K_{n_{d-1}}$.
Therefore, by induction on $d$, we assume that we can put $2k T_{d-1}$
pebbles on $(x_1, x_2, \ldots, x_{d-1}, v)$.  But now we have $2k
T_{d-1}$ pebbles on $(x_1, x_2, \ldots, x_{d-1}, x_d)$ and we can move
an additional $k T_{d-1}$ pebbles from $(x_1, x_2, \ldots, x_{d-1},
v)$ to $(x_1, x_2, \ldots, x_{d-1}, x_d)$ for every vertex $v \neq
x_d$.  Thus, we can move a total of $(n_d+1) k T_{d-1} = k T_d$
pebbles onto $(x_1, x_2, \ldots, x_{d-1}, x_d)$, as required, and so
$\pi_{kT_d}^* (G) \leq 2^d k \displaystyle{\prod_{i=1}^d n_i}$.

Conversely, we know from Theorem~\ref{bounds for pi^*(GxKn)} that
\[
\pi_{kT_d}^* (K_{n_1} \cart K_{n_2} \cart \cdots \cart K_{n_d}) \geq
\left\lceil \left( \frac{n_d}{n_d+1} \right) \pi^*_{2k T_d}(K_{n_1}
\cart K_{n_2} \cart \cdots \cart K_{n_{d-1}}) \right \rceil.
\]
Now $2k T_d = 2k (n_d+1) T_{d-1}$, and we may assume by induction on
$d$ that
\[ \pi^*_{2k (n_d+1)T_{d-1}}(K_{n_1} \cart K_{n_2} \cart \cdots \cart
K_{n_{d-1}}) = 2^{d-1} (2k (n_d+1)) \displaystyle{\prod_{i=1}^{d-1}
n_i} = 2^d k (n_d+1) \displaystyle{\prod_{i=1}^{d-1} n_i}.
\]
Multiplying this number by $\frac{n_d}{n_d+1}$ gives an integer, so
taking the ceiling is irrelevant.  Therefore,
\[
\pi_{kT_d}^* (K_{n_1} \cart K_{n_2} \cart \cdots \cart K_{n_d}) \geq
2^d k n_d \displaystyle{\prod_{i=1}^{d-1} n_i} = 2^d k
\displaystyle{\prod_{i=1}^d n_i},
\]
which agrees with our upper bound.  Therefore, $\pi_{kT_d}^* (K_{n_1}
\cart K_{n_2} \cart \cdots \cart K_{n_d}) = 2^d k
\displaystyle{\prod_{i=1}^d n_i}$.\Bx
\begin{theorem}
\label{asymptotic bound for complete products}
Let $n_1, n_2, \ldots, n_d$ be a sequence of nonnegative integers.
Then
\[
\pi_t^* (K_{n_1} \cart K_{n_2} \cart \cdots \cart K_{n_d}) \in 2^d t
\prod_{i=1}^d \frac{n_i}{n_i + 1} + \Theta(1).
\]
\end{theorem}
\textbf{Proof}: Let $G = K_{n_1} \cart K_{n_2} \cart \cdots \cart
K_{n_d}$, let $T = \displaystyle{\prod_{i=1}^d (n_i + 1)}$, and let
$C$ be given by
\[
C = \max_{t < T} \left( \pi_t^* (G) - 2^d t \prod_{i=1}^d
\frac{n_i}{n_i + 1} \right).
\]
For any $t > T$, we can let $q = t \mbox{ div } T$ and $r = t
\bmod{T}$.  Then applying Proposition~\ref{pi^*_s+t}, we have
\[
\pi_t^*(G) = \pi_{qT+r}^*(G) \leq \pi_{qT}^*(G) + \pi_r^*(G) \leq q
\pi_T^*(G) + \pi_r(G).
\]
From Lemma~\ref{2^d pebbles on each vertex}, we know that
$\pi_{qT}^*(G) = 2^d q \displaystyle{\prod_{i=1}^d n_i} = 2^d q T
\displaystyle{\prod_{i=1}^d \frac{n_i}{n_i + 1}}$, and from the
definition of $C$, we have
\[
\pi_t^*(G) \leq 2^d q T \prod_{i=1}^d \frac{n_i}{n_i+1} +
2^d r \prod_{i=1}^d \frac{n_i}{n_i + 1} + C = 2^d (qT+r) \prod_{i=1}^d
\frac{n_i}{n_i + 1} + C = 2^d t \prod_{i=1}^d \frac{n_i}{n_i + 1} + C.
\]
Thus, $0 \leq \pi_t^* (G) - 2^d t \displaystyle{\prod_{i=1}^d
\frac{n_i}{n_i + 1}} \leq C$, and so $\displaystyle{\pi_t^* (G) \in
2^d t \prod_{i=1}^d \frac{n_i}{n_i + 1} + \Theta(1)}$, as desired.\Bx

\section{Optimal Fractional Pebbling}
\label{fractional}

\emph{Fractional distributions} and \emph{fractional pebbling moves}
were defined in~\cite{diameter fractional and target}.  These are
continuous analogs of pebbling concepts.  Moews~\cite{optimal
hypercubes} previously called them \emph{continuous distributions},
and \emph{continuous pebbling moves}, and he defined the
\emph{continuous optimal pebbling number} of a graph.  We give these
definitions now.

\noindent
\textbf{Definitions \cite{diameter fractional and target}}: A
\emph{fractional distribution} on $G$ is a function $D: V \rightarrow
\mathbb{R}^+ \cup \{ 0 \}$.  Again, a distribution represents a
placement of pebbles on the vertices of $G$, though we now allow a
nonintegral number of pebbles.  A \emph{fractional pebbling move}
consists of removing $2k$ pebbles from one vertex and adding $k$
pebbles to an adjacent vertex.  As in an integer-valued distribution,
the \emph{size} of $D$ is given by $|D| = \displaystyle{\sum_{v\in V}
D(v)}$, and $D$ is \emph{fractionally solvable}, (or simply
\emph{solvable} if there is no ambiguity), in case for every vertex
$v$, it is possible to reach $v$ with one pebble through some sequence
of fractional pebbling moves, starting from $D$.

Moews~\cite{optimal hypercubes} defined the \emph{continuous optimal
pebbling number} of a graph, which we denote $\overline{\pi}^* (G)$.
The \emph{optimal fractional pebbling number} of the graph $G$, which
we denote $\hat{\pi}^*(G)$, was defined in~\cite{diameter fractional
and target}.  We give these definitions now.

\noindent
\textbf{Definitions \cite{diameter fractional and target, optimal
hypercubes}}: The \emph{continuous optimal pebbling number} of a graph
$G$, is the smallest number $\overline{\pi}^* (G)$ such that some
fractional distribution $D$ with $|D| = \overline{\pi}^*(G)$ is
solvable using fractional pebbling moves.  The \emph{optimal
fractional pebbling number} $\hat{\pi}^*(G)$ is given by
\[
\hat{\pi}^*(G) = \liminf_{t \rightarrow \infty}
\frac{\pi_{t}^*(G)}{t}.
\]
Theorem~\ref{optimal fractional theorem} was shown in~\cite{diameter fractional and target}.
\begin{theorem}[\cite{diameter fractional and target}]
\label{optimal fractional theorem}
Every graph $G$ satisfies $\hat{\pi}^*(G)=\overline{\pi}^*(G)$.  Furthermore, $\hat{\pi}^*(G)$ is rational for any graph $G$ and every graph has an optimal fractional distribution in which the number of pebbles on each vertex is rational.\Bx
\end{theorem}
%It was shown in~\cite{diameter fractional and target} that
%$\hat{\pi}^*(G) = \overline{\pi}^*(G)$ for every graph $G$.
%Furthermore, as with invariants like fractional chromatic number,
%fractional poset dimension, fractional matching number, and others,
%$\hat{\pi}^*(G)$ is equal to the solution of a particular linear
%optimization problem with integer coefficients, and thus is rational
%for any graph.  Also, there is an optimal fractional distribution in 
%which the number of pebbles on each vertex is rational.
%
Moews~\cite{optimal hypercubes} proved Theorem~\ref{continuous optimal graham}, and used
it to give a nonconstructive proof of Theorem~\ref{optimal products},
which relates the optimal pebbling number of $G^d$ to the continuous
optimal pebbling number of $G$.
\begin{theorem}[Moews~\cite{optimal hypercubes}]
\label{continuous optimal graham}
For all graphs $G$ and $G'$, we have $\hat{\pi}^* (G \cart G') =
\hat{\pi}^* (G) \hat{\pi}^* (G')$.
\end{theorem}
\begin{theorem}[Moews~\cite{optimal hypercubes}]
\label{optimal products}
For all graphs $G$, we have $\pi^* (G^d) \in
O(\hat{\pi}^*(G)^{d+c \log d}) = O(\left( \hat{\pi}^*(G)
\right)^d \cdot d^k)$ for some constants $c$ and $k$.
\end{theorem}

Theorem \ref{optimal fractional} generalizes some of our results from
Section~\ref{kn products section}.  We begin with Lemma \ref{optimal
fractional lemma}.

\begin{lemma}\label{optimal fractional lemma}
For every graph $G$ we have $\pi_t^*(G)\geq \hat{\pi}^*(G)t$ for all $t$.
\end{lemma}
\textbf{Proof}: Suppose by contradiction that there is some $t$ such
that $\pi_t^*(G)< \hat{\pi}^*(G)t$.  Let $D$ be a $t$-solvable
distribution on $G$ with $|D|=\pi_t^*(G)$.  Then the fractional
distribution $\hat{D}$ given by $\hat{D}(v)=\frac{D(v)}{t}$ for all
$v$ is fractionally solvable.  So,
$\overline{\pi}^*(G)\leq\frac{\pi_t^*(G)}{t}<\hat{\pi}^*(G)$,
contradicting Theorem~\ref{optimal fractional theorem}. \Bx \\
The definition of $\hat{\pi}^*(G)$ implies that
$\pi_t^*(G) \in \hat{\pi}^*(G)t + o(t)$.  Here we tighten the lower
order term.
\begin{theorem}
For every graph $G$ we have $\pi_t^*(G) \in \hat{\pi}^*(G)t+\Theta(1)$.
\label{optimal fractional}
\end{theorem}
\textbf{Proof}: Let $V(G)=\{v_1,v_2,\ldots,v_n\}$ and let $a$ and $b$
be integers satisfying $\hat{\pi}^*(G)=\frac{a}{b}$.  From Theorem~\ref{optimal fractional theorem}, there is some fractionally
solvable fractional distribution $\hat{D}$ on $G$ with $|\hat{D}| =
\frac{a}{b}$ such that $\hat{D}(v_i)=\frac{a_i}{b_i}$ for some integers $a_i$ and $b_i$.
Under $\hat{D}$, if ${\rm dist}(v_i,v_j)=\delta$, then $\left(
\frac{a_i}{b_i} \right) 2^{-\delta}$ pebbles could be sent from $v_i$
to $v_j$ by making fractional pebbling moves toward $v_j$.  Since
$\hat{D}$ is fractionally solvable, every vertex $v$ satisfies
\[ \sum_i \left(\frac{a_i}{b_i}\right) 2^{-{\rm dist}(v_i,v)} \ge 1. \]
Let $l={\rm lcm}(b_1,b_2,\ldots,b_n)$ and let $k=2^d l$, where $d={\rm
diam}(G)$.  Given an integer $t$, the division algorithm produces
integers $q$ and $r$ such that $t=kq+r=2^d lq+r$ and $0\leq r\leq
k-1$.  Consider the distribution $D$ on $G$ given by $D(v_i) = kq
\hat{D}(v_i)+r$ for all $i$.  Under $D$, we have $D(v_i) = \left(
\frac{a_i}{b_i} \right) 2^d lq+r$ for all $i$.  Since $l$ is a
multiple of $b_i$, $D'(v_i)=D(v_i)-r$ is a multiple of $2^d$.  So,
under the distribution $D'$, it is possible to send $\left(
\frac{a_i}{b_i} \right) 2^{d-\delta} lq$ pebbles from $v_i$ to $v_j$.
So, starting from $D'$, the number of pebbles that can be sent to a
root $v$ is given by
\[
\sum_i \left( \frac{a_i}{b_i} \right) 2^{d-{\rm dist}(v_i,v)} lq =
2^d lq \sum_i \left( \frac{a_i}{b_i} \right) 2^{-{\rm dist}(v_i,v)} \ge
2^d lq = t-r.
\]
Thus, $D'$ is $(t-r)$-solvable on $G$, meaning $D$ is
$t$-solvable on $G$.  Since $n$ and $k$ are constants, we have
\[
\pi_t^*(G) \ \leq \ kq\left(\frac{a}{b}\right)+nr \ \ \leq \
\left(\frac{a}{b}\right)t+n(k-1) \in \
\hat{\pi}^*(G)t+O(1).
\]
In connection with Lemma \ref{optimal fractional lemma}, this gives us
the desired result. \Bx \\
We note that $\hat{\pi}^* (K_n) = \frac{2n}{n+1}$, $\hat{\pi}^* (K_2
\cart K_2)=\frac{16}{9}$, and $\hat{\pi}^*(K_2 \cart K_3) = 2 =
\frac{12}{6}$.  Thus, these specific cases of Theorem \ref{optimal
fractional} are witnessed by Theorem \ref{pi t(Kn)} and Propositions
\ref{pi t(K2xK2)} and \ref{pi t(K2xK3)}, respectively.

\section{Products of $C_5$}
\label{C5 x ... x C5}

If we apply Theorem~\ref{optimal products} to $C_5$, we find
$\pi^*(C_5^d) \in O(2^d d^k)$ for some constant $k$, since
$\hat{\pi}^*(C_5) = 2$.  However, Moews's proof of
Theorem~\ref{optimal products} was nonconstructive.  It does not give
distributions for small values of $d$, and it gives no information for
small values of $d$.  We give distributions that show that
$\pi^*(C_5^d) \in O(\sqrt{5}^d)$.  We let the vertices of $C_5$ be $\{
v_0, v_1, v_2, v_3, v_4 \}$.  We begin by finding $t$-solvable
distributions $A_t$ on $C_5 \cart C_5$ for $t=1$, $t=2$, and $t=4$.

\begin{figure}
\putpicture{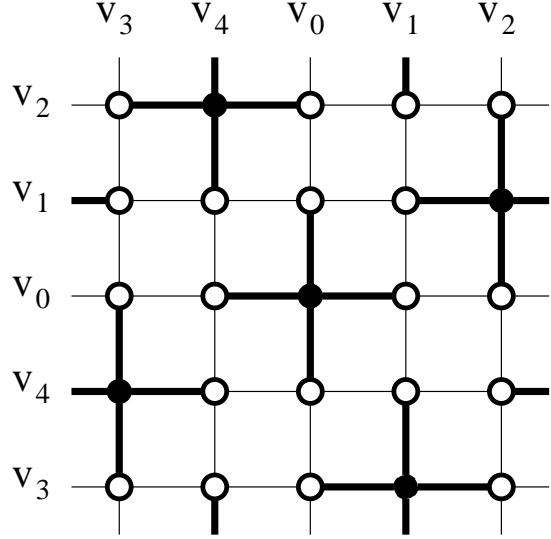}{2.8in}
\caption{The distribution $B$ on $C_5 \cart C_5$}
\label{B}
\end{figure}

\noindent
\textbf{Notation}: We denote by $A_1$ the distribution with four
pebbles on $(v_0, v_0)$ and two pebbles each on $(v_2, v_2)$ and
$(v_3, v_3)$, by $A_2$ the distribution with four pebbles each on
$(v_0, v_0)$, $(v_2, v_2)$, and $(v_3, v_3)$, and by $A_4$ the
distribution with 4 pebbles on each $(v_i,v_{2i\bmod{5}})$, $0\le i\le
4$.  We write $B$ for the $\frac{1}{4}A_4$, i.~e.\ the distribution
with one pebbles on each $(v_i, v_{2i\bmod{5}})$.  $B$ is shown in
Figure~\ref{B} (filled-in vertices are occupied, dark edges give the
neighborhoods of the occupied vertices, and the edges wrap around in
the obvious ways).
\begin{proposition}
\label{t-solvable on C5xC5, t=1,2,4}
For each $t\in\{1,2,4\}$, the distribution $A_t$ is $t$-solvable on 
$C_5 \cart C_5$.
\end{proposition}
\textbf{Proof}: First note that $\{ v_2, v_3 \} \cart \{ v_2, v_3 \}
\cong K_2 \cart K_2$.  We call these vertices the \emph{corners} of
the graph, imagining $(v_0, v_0)$ to be the center.  If we have two
pebbles each on $(v_2, v_2)$ and $(v_3, v_3)$, we have the
$2$-solvable distribution in Proposition~\ref{pi t(K2xK2)}, so two
pebbles can be moved any of the corners, and one pebble can be moved
to any vertex adjacent to these corners.  The rest of the vertices are
within two steps from $(v_0, v_0)$, so they can be reached from the
four pebbles from there.  This takes care of the $t=1$ case.

For $t=2$, we instead have four pebbles each on $(v_2, v_2)$ and
$(v_3, v_3)$, so we can consider these to be two groups which each have
two pebbles on both $(v_2, v_2)$ and $(v_3, v_3)$.  Therefore, we can
put four pebbles on any corner.  Then the vertices whose distance from
$(v_0, v_0)$ is zero or one can receive two pebbles from that vertex.
The vertices whose distance from $(v_0, v_0)$ is three or four can
receive two pebbles from the nearest corner, and those whose distance
from $(v_0, v_0)$ is two can receive one pebble from $(v_0, v_0)$ and
one from the nearest corner.

When $t=4$, the symmetry of $A_4$ allows us to consider only one
target, say $(v_0,v_1)$.  This vertex can receive two pebbles from
$(v_0, v_0)$, and one each from $(v_3, v_1)$ and $(v_1, v_2)$.  \Bx
\begin{theorem}
\label{B dot D}
Let $G$ be any graph and let $D$ be a $t$-solvable distribution on $G$
in which the number of pebbles on every vertex is a multiple of four.
Then the distribution $B \cdot D$ is a $t$-solvable distribution in
$(C_5 \cart C_5) \cart G$ in which the number of pebbles on every
vertex is a multiple of four.  Note that the number of pebbles in $B
\cdot D$ is $5 |D|$.  In particular, by induction on $m$, we have
$\pi_t^*(C_5^{2m} \cart G) \leq 5^m |D|$.
\end{theorem}
\textbf{Proof}: Let the target vertex in $C_5 \cart C_5 \cart
G$ be $(v_i, v_j, y)$.  Since $D(x)$ is a multiple of four, we write
$\frac{1}{4} D$ for the distribution with $\frac{1}{4} D(x)$ pebbles
on $x \in V$, and we write $\mathcal{S}_{1/4}$ for the set
$\mathcal{S}_{1/4} = \left \{ \frac{1}{4} D \right \}$.  Now
$B \cdot D = (4B) \cdot \left( \frac{1}{4} D \right)$.  By
Proposition~\ref{t-solvable on C5xC5, t=1,2,4},
$4B$ is $\mathcal{S}_4(C_5 \cart
C_5)$-solvable, so by Theorem~\ref{optimal pebbling distributions}, $B
\cdot D$ is $\left( \mathcal{S}_4(C_5 \cart C_5) \cdot
\mathcal{S}_{1/4} \right)$-solvable in $C_5 \cart C_5 \cart
G$.  That is, from $B \cdot D$ we can reach the distribution in which
$4 \left( \frac{1}{4} D(y_k) \right) = D(y_k)$ pebbles are on $(v_i,
v_j, y_k)$ for every $y_k \in V$.  But now the distribution on the
vertices in $(v_i, v_j) \cart G \cong G$ is $D$.  Since $D$ is
$t$-solvable, we can put $t$ pebbles on $(v_i, v_j, y_k)$.  Clearly,
$(B\cdot D) ((v_i, v_j, y_k)) = B((v_i, v_j)) D(y_k)$ is a multiple of
four, since $D(y_k)$ is a multiple of four.\Bx 
%
%\begin{corollary}
%\label{odd product of C5}
%For all integers $m \geq 0$, we have:
%\begin{eqnarray*}
%\pi^*(C_5^{2m+1}) \leq 4 \cdot 5^m \\
%\pi_2^*(C_5^{2m}) \leq \frac{12}{5} \cdot 5^m \\
%\pi_4^* (C_5^{2m}) \leq 4 \cdot 5^m.
%\end{eqnarray*}
%\end{corollary}
\begin{corollary}
\label{odd product of C5}
For all integers $m \geq 0$, we have $\pi^*(C_5^{2m+1}) \leq 4 \cdot 5^m$.
\end{corollary}
%
%(Of course, one can also obtain $\pi_2^*(C_5^{2m}) \leq \frac{12}{5}
%\cdot 5^m$ and $\pi_4^* (C_5^{2m}) \leq 4 \cdot 5^m$, but we will not
%need these results.)\\
%
\textbf{Proof}: We apply Theorem~\ref{B dot D} to the distribution
with four pebbles on a single vertex of $G = C_5$.
%For $\pi_2^*(C_5^{2m})$ with $m \geq 1$, apply Theorem~\ref{B dot D}
%to the distribution $A_2$.  For $\pi_4^* (C_5^{2m})$, apply
%Theorem~\ref{B dot D} to the trivial graph with four pebbles on the
%lone vertex. 
\Bx

A natural question at this point is what bounds we can get for
$\pi^*(C_5^{2m})$.  We create a solvable distribution $F$ on $C_5^4$,
and we use $F$ to start an induction with Theorem~\ref{B dot D} for
even products similar to the argument for Corollary~\ref{odd product
of C5}. 

\noindent
\textbf{Notation}: Let $F$ be the distribution of $44$ pebbles on
$C_5^4$ given by
\[ F(v_i, v_j, v_k, v_m) = \left \{ \begin{array}{cl}
A_4(v_k, v_m) & \mbox{if } i=j=0 \\
A_2(v_k, v_m) & \mbox{if } i=j=2 \mbox{ or } i=j=3 \\
0 & \mbox{otherwise.}
\end{array} \right. \]
Note that if we denote the empty distribution by $A_0$, then we may
more simply write $$F(v_i,v_j,v_k,v_m)=A_r(v_k,v_m),$$ 
where $r=A_1(v_i,v_j)$.
\begin{proposition}
\label{F is 4-solvable}
Every occupied vertex in $F$ has four pebbles, and $F$ is solvable in
$C_5^4$.
\end{proposition}
\textbf{Proof}: Every occupied vertex of both $B$ and $A_2$ has four
pebbles in $C_5 \cart C_5$, so this holds for $F$ in $C_5^4$ as well.
To show $F$ is solvable, let the target vertex in $C_5^4$ be $(v_i,
v_j, v_k, v_m)$.  By construction, the distribution of pebbles on
$(v_0, v_0) \cart C_5 \cart C_5 \cong C_5 \cart C_5$ is $B$.
Therefore, by Proposition~\ref{t-solvable on C5xC5, t=1,2,4},
four pebbles can be
moved to $(v_0, v_0, v_k, v_m)$ using only the pebbles on $(v_0, v_0)
\cart C_5 \cart C_5$.  Similarly, and simultaneously, by
Proposition~\ref{t-solvable on C5xC5, t=1,2,4}, two pebbles can be
moved to both $(v_2, v_2, v_k, v_m)$ and $(v_3, v_3, v_k, v_m)$.  At
this point, the distribution of pebbles on $C_5 \cart C_5 \cart
(v_k, v_m) \cong C_5 \cart C_5$ is $A_1$, so one pebble may be moved
to $(v_i, v_j, v_k, v_m)$, again by Proposition~\ref{t-solvable on
C5xC5, t=1,2,4}.\Bx

The use of Theorem~\ref{B dot D} on the distribution $A_1$ would give
a better coefficient of $\frac{8}{5}$ in the Theorem~\ref{all products
of C5}; however, $A_1$ does not qualify since some vertices get only 2
pebbles.
\begin{theorem}
\label{all products of C5}
We have $\pi^*(C_5^{2m}) \leq \frac{44}{25} (5^m)$.  For all $d \geq
1$ we have $\pi^*(C_5^d) \leq \frac{4}{\sqrt{5}} (5^{\frac{d}{2}}) \in
O(\sqrt{5}^d)$.
\end{theorem}
\textbf{Proof}: If $d=2$, Proposition~\ref{t-solvable on C5xC5, t=1,2,4} 
shows that $\pi^* (C_5 \cart C_5) \leq 8 < \frac{44}{5}$.  If
$d=4$, Proposition~\ref{F is 4-solvable} shows that $\pi^* (C_5^4)
\leq 44$.  For $d=2m$ with $m > 2$, Theorem~\ref{B dot D} implies that
$\pi^* (C_5^{2m}) = \pi^* (C_5^{2(m-2)} \cart C_5^4) \leq 5^{m-2} |F|
= \frac{44}{25} (5^m)$.  Since $\frac{44}{25} < \frac{4}{\sqrt{5}}$,
the second part follows for even $d$.

If $d=2m+1$, Corollary~\ref{odd product of C5} gives us $\pi^*(C_5^d)
\leq 4 \cdot 5^m = 4 \cdot 5^{\frac{n-1}{2}} = \frac{4}{\sqrt{5}}
(5^{\frac{n}{2}})$.  \Bx 

We can generalize the construction of $F$ and the proof of
Proposition~\ref{F is 4-solvable} to obtain Theorem~\ref{using a
solvable distribution}.
\begin{theorem}
\label{using a solvable distribution}
Let $\mathcal{S}$ be a set of distributions on $G$, suppose $D$ is an
$\mathcal{S}$-solvable distribution, and suppose $\{ D_r' \}_{r \geq
1}$ is a family of distributions on $G'$ such that each $D'_r$ is
$r$-solvable.  Let $\Delta : V(G \cart G') \rightarrow \mathbb{N}$
be the distribution on $G \cart G'$ defined by
\[ \Delta((v, w)) = D'_{D(v)} (w). \] 
Then $\Delta$ is $(\mathcal{S} \cdot \mathcal{S}_1 (G'))$-solvable
in $G \cart G'$.  That is, for any distribution $\overline{D} \in
\mathcal{S}$ and any vertex $w \in V'$, a copy of $\overline{D}$ can
be moved to the vertices of $G \cart \{ w \}$.
\end{theorem}
\textbf{Proof}: Let $w$ be the chosen vertex in $G'$.  Then for any $v
\in V$, restricting $\Delta$ to the vertices $\{ v \} \cart G'$ gives
the distribution $D_{D(v)}$ in $\{ v \} \cart G'$.  Since $D_{D(v)}$
is $D(v)$-solvable in $G'$, we can move $D(v)$ pebbles to $(v, w)$ for
each $v \in V$.  After these moves, the distribution of pebbles on $G
\cart \{ w \}$ is $D$.  Since $D$ is $\mathcal{S}$-solvable, we can
put a copy of any $\overline{D} \in \mathcal{S}$ on $G \cart \{ w
\}$, as desired.\Bx 

We note that the proofs of Corollaries~\ref{optimal st-pebbling
graphs} and~\ref{optimal pebbling graphs} essentially involved letting
each $D_r = r D'$, where $D'$ is the distribution of pebbles on
$G'$ defined in the proof of Theorem~\ref{optimal pebbling
distributions}.  Corollary~\ref{optimally using solvable distribution}
is a stronger result.
\begin{corollary}
\label{optimally using solvable distribution}
Let $\mathcal{S}$ be a set of distributions on $G$, and suppose
$D$ is an $\mathcal{S}$-solvable distribution.  Then for any graph
$G'$, we have
\[
\pi^* (G \cart G', \mathcal{S} \cdot \mathcal{S}_1(G')) \leq \sum_{v
\in V} \pi_{D(v)}^* (G').
\]
\end{corollary}
\textbf{Proof}: We simply apply Theorem~\ref{using a solvable
distribution} and use a family of distributions $\{ D_t' \}$ in
which each $D_t'$ is optimal, i.~e.\ $|D_t'| = \pi_t^*(G')$.\Bx

\section{Hypercubes}
\label{hypercubes section}

In this section, we give optimal pebbling distributions on the
$d$-dimensional hypercube $Q^d \cong K_2^d$.  We consider the vertices
of $Q^d$ to be all bitstrings of length $d$, or equivalently, all
vectors in the $d$-dimensional vector space $\mathbb{F}_2^d$ over the
two-element field $\mathbb{F}_2$.  There is an edge between two
vertices when the Hamming distance between the corresponding
bitstrings is $1$.  Given two bitstrings $\mathbf{v_1} \in V(Q^{d_1})$
and $\mathbf{v_2} \in V(Q^{d_2})$, we write $\mathbf{v_1} \cdot
\mathbf{v_2}$ for the bitstring in $V(Q^{d_1+d_2})$ obtained by
concatenating the bits in $\mathbf{v_1}$ and $\mathbf{v_2}$.  We also
write $\mathbf{0^k}$ and $\mathbf{1^k}$ for the bitstrings $00 \ldots
0$ and $11 \ldots 1$, respectively, and we call the number of $1$'s in
a bitstring its \emph{weight}.

Since the continuous optimal pebbling number of $K_2$ is
$\hat{\pi}^* (K_2) = \frac{4}{3}$, Theorem~\ref{optimal products}
implies Theorem~\ref{Moews bound}, which Moews also proved directly.
\begin{theorem}[Moews~\cite{optimal hypercubes}]
\label{Moews bound}
The optimal pebbling number of $Q^d$ satisfies $\pi^* (Q^d) \in
O \left( \frac{4}{3}^d d^k \right)$ for some constant $k$.
\end{theorem}
Theorem~\ref{Moews bound} gives the best known bound for hypercubes,
but it does not give explicit distributions, which is our aim.  The
$d^{th}$ root of Moews's result tends to about $1.33$, and
Proposition~\ref{antipodal solvable} gives an example, the $d^{th}$
root of whose size is roughly $1.41$.  Our new construction in
Theorem~\ref{weighted geometric mean} improves that number below
$1.38$.  Proposition~\ref{antipodal solvable} gives a solvable
distribution on $Q^d$ for all $d$.  These were first given in Pachter,
Snevily, and Voxman~\cite{PSV}.
\begin{proposition}[Pachter et al.~\cite{PSV}]
\label{antipodal solvable}
If $d = 2k$, the distribution on $Q^d$ obtained by putting $2^k$
pebbles on $\mathbf{0^d}$ and $2^{k-1}$ pebbles on $\mathbf{1^d}$ is
solvable.  If $d = 2k+1$, the distribution on $Q^d$ given by putting
$2^k$ pebbles on both $\mathbf{0^d}$ and $\mathbf{1^d}$ is solvable.
Thus, the optimal pebbling number of a hypercube satisfies
\[ \begin{array}{c}
\pi^*(Q^{2k}) \leq 3 \cdot 2^{k-1} \\
\pi^*(Q^{2k+1}) \leq 2^{k+1}. \\
\end{array} \]
In particular, $\pi^*(Q^d) \in O(2^{\frac{d}{2}}) = O(\sqrt{2}^d)$.
\end{proposition}
\textbf{Proof}: In both cases whether $d$ is even or odd, every vertex
whose weight is at most $k$ can receive at least one pebble from the
pebbles on $\mathbf{0^k}$ in the given distribution, and every vertex
with larger weight can receive a pebble from those on
$\mathbf{1^k}$.\Bx

We give a construction for extending the distributions in
Proposition~\ref{antipodal solvable} to distributions on larger cubes
with better asymptotic bounds than those in the Proposition.
This construction is based on an argument similar to the proof
of Theorem~\ref{using a solvable distribution} using distributions on
$K_2$ obtained from Theorem~\ref{pi t(Kn)}.  First recall the
distributions on $K_2$ from Theorem~\ref{pi t(Kn)}; we will use these
in Theorem~\ref{splitting construction}.

\noindent
\textbf{Definition}: We let $\mathcal{D}$ be the family $\mathcal{D} =
\{ D_r \}_{r \geq 1}$ of distributions on $K_2$ given by
\[ \begin{array}{c*{2}{@{\hspace{1in}}c}}
D_{3k}(x_0) = 2k & D_{3k+1}(x_0) = 2k+2 & D_{3k+2}(x_0) = 2k+2 \\
D_{3k}(x_1) = 2k & D_{3k+1}(x_1) = 2k & D_{3k+2}(x_1) = 2k+1 \\
\end{array} \]
\begin{proposition}
\label{Dr r-solvable}
Each $D_r$ is $r$-solvable, and in each case, we have $|D_r| = \left
\lceil \frac{4r}{3} \right \rceil \leq \frac{4}{3}r+\frac{2}{3}$.
\end{proposition}
\textbf{Proof}: Each $D_r$ is the $r$-solvable distribution from the
proof of Theorem~\ref{pi t(Kn)}.  Counting pebbles, we have $|D_{3k}|
= 4k$, $|D_{3k+1}| = 4k+2$, and $|D_{3k+2}| = 4k+3$.  In each case,
$|D_r| = \left \lceil \frac{4}{3} r \right \rceil$.\Bx

In the spirit of Theorem~\ref{using a solvable distribution}, we want
to extend a solvable distribution $D$ on a graph $G$ to a distribution
$D'$ on $G \cart K_2$.  We hope that $|D'| \approx \frac{4}{3} |D|$.
Unfortunately, the extra $\frac{2}{3}$ in Proposition~\ref{Dr
r-solvable} can cause problems.  For example, if $D$ has a single
pebble on a large number of vertices, those pebbles each give rise to
two pebbles in $D'$.  We can get an extra $\frac{2}{3}$ for each
occupied vertex in $D$.  We define the \emph{support} of $D$ to keep
track of this information. 

\noindent
\textbf{Definition}: The \emph{support} of a distribution $D$ on the
graph $G$, denoted $\sigma(D)$, is the set of occupied vertices in
$D$; i.~e.\ $\sigma(D) = \{ v \in V(G) : D(v) > 0 \}$.
\begin{theorem}
\label{splitting construction}
Let $D$ be a $t$-solvable distribution on $Q^d$.  For each $\mathbf{v}
\in V(Q^d)$, define $D'(\mathbf{v} \cdot 0)$ and $D'(\mathbf{v} \cdot
1)$ by
\[ \begin{array}{c}
D'(\mathbf{v} \cdot 0) = D_{D(\mathbf{v})}(x_0) \\
D'(\mathbf{v} \cdot 1) = D_{D(\mathbf{v})}(x_1) \\
\end{array} \]
Then $D'$ is $t$-solvable on $Q^{d+1}$.  Furthermore, the number of
pebbles in $D'$ is at most $\frac{4}{3} |D| + \frac{2}{3}
|\sigma(D)|$, and $|\sigma(D')| \leq 2 |\sigma(D)|$.
\end{theorem}
\textbf{Proof}: Let the target in $Q^{d+1}$ be $\mathbf{v} \cdot b$,
where $\mathbf{v} \in V(Q^d)$ and $b \in \{ 0, 1 \}$.  For each
$\mathbf{v_i} \in V(Q^d)$, the distribution of pebbles on
$\mathbf{v_i} \cart K_2 \cong K_2$ is $D_{D(\mathbf{v_i})}$.  Since
this distribution is $D(\mathbf{v_i})$-solvable in $K_2$, we can put
$D(\mathbf{v_i})$ pebbles on $\mathbf{v_i} \cdot b$.  If we do this
for each $\mathbf{v_i} \in V(Q^d)$, the distribution of pebbles on
$Q^d \cart \{ b \} \cong Q^d$ is $D$.  Since $D$ is $t$-solvable on
$Q^d$, we can put $t$ pebbles on $\mathbf{v} \cdot b$.  The total
number of pebbles in $D'$ is
\[
|D'| = \sum_{\mathbf{v}\in \sigma(D')} D'(\mathbf{v}) =
\sum_{\mathbf{v}\in \sigma(D)} |D_{D(\mathbf{v})}| \leq
\sum_{\mathbf{v}\in \sigma(D)} \left( \frac{4}{3} D(\mathbf{v}) +
\frac{2}{3} \right) = \frac{4}{3} |D| + \frac{2}{3} |\sigma(D)|.
\]
Finally, $\sigma(D') \subseteq \sigma(D) \cart \{ 0, 1 \}$.\Bx

Theorem~\ref{many splits} describes what happens when we apply
Theorem~\ref{splitting construction} repeatedly.
\begin{theorem}
\label{many splits}
Let $D$ be a solvable distribution on a graph $G$ with $s =
|\sigma(D)|$, and let $D_m$ be the result of applying
Theorem~\ref{splitting construction} $m$ times to $D$.  Then $D_m$ is
a solvable distribution on $G \cart K_2^m \cong G \cart Q^m$ such that
$|\sigma (D_m)| \leq 2^m s$, and $|D_m| \leq \left( \frac{4}{3}
\right)^m |D| + 2^m s - \left( \frac{4}{3} \right)^m s$.
\end{theorem}
\textbf{Proof}: There is nothing to show if $m=0$, so we suppose by
induction that for some $i \geq 0$, $D_i$ is a solvable distribution
on $G \cart Q^i$ with $|D_i| \leq \left( \frac{4}{3} \right)^i |D| +
2^i s - \left( \frac{4}{3} \right)^i s$, and $|\sigma (D_i)| \leq 2^i
|\sigma (D)|$.  Then applying Theorem~\ref{splitting construction} to
$D_i$, we find that $D_{i+1}$ is a solvable distribution on $G \cart
Q^{i+1}$ with $|\sigma (D_{i+1})| \leq 2 |\sigma(D_i)| = 2^{i+1} s$.
Furthermore, we have
\[
|D_{i+1}| \leq \frac{4}{3} |D_i| + \frac{2}{3} |\sigma(D_i)| \leq
\frac{4}{3} \left[ \left( \frac{4}{3} \right)^i |D| + 2^i s - \left(
\frac{4}{3} \right)^i s \right] + \frac{2}{3} (2^i s).
\]
Multiplying through by the $\frac{4}{3}$ and noting that $\frac{4}{3}
(2^i s) + \frac{2}{3} (2^i s) = 2^{i+1}s$, we have
\[
|D_{i+1}| \leq \left( \frac{4}{3} \right)^{i+1} |D| + \frac{4}{3} (2^i
s) - \left( \frac{4}{3} \right)^{i+1} s + \frac{2}{3} (2^i s)
= \left( \frac{4}{3} \right)^{i+1} |D| + 2^{i+1} s - \left(
\frac{4}{3} \right)^{i+1} s,
\]
completing the induction.\Bx
\begin{corollary}
\label{splitting antipodal}
Let $D_0$ be the distribution with $2^k$ pebbles on both
$\mathbf{0^{2k+1}}$ and $\mathbf{1^{2k+1}}$ in $Q^{2k+1}$, and let
$D_m$ be the resulting distribution in $Q^{2k+m+1}$ obtained by
applying Theorem~\ref{splitting construction} $m$ times.  Then
\[
|D_m| \leq \left( \frac{4}{3} \right)^m (2^{k+1}) + 2^{m+1} - 2\left(
 \frac{4}{3} \right)^m.
\]
\end{corollary}
\textbf{Proof}: We apply Theorem~\ref{many splits} to $D_0$, noting
that $|D_0| = 2^{k+1}$ and $s = |\sigma(D_0)| = 2$.\Bx 

For large $m$, the term $2\left( \frac{4}{3} \right)^m$ in
Corollary~\ref{splitting antipodal} is small compared to $2^{m+1}$.
By controlling the relationship between $k$ and $m$, we can ensure
that the first two terms are roughly equal.  Using logarithms to solve
the equation $\left( \frac{4}{3} \right)^m (2^{k+1}) \approx 2^{m+1}$,
or $2^k \approx \left( \frac{3}{2} \right)^m = 1.5^m$ together with
the observation that $d = 2k+m+1$. we obtain the constants in
Theorem~\ref{weighted geometric mean}.
\begin{theorem}
\label{weighted geometric mean}
Given an integer $d$, let $k = \left \lceil \frac{\log_2 1.5}{\log_2
  4.5} (d-1) \right \rceil \approx 0.2696(d-1)$, and let $m = d - 1 - 2k
\approx 0.4608 (d-1)$.  Then the distribution $D_m$ on $Q^{2k+m+1} =
Q^d$ from Corollary~\ref{splitting antipodal} with these values of $k$
and $m$ satisfies $|D_m| \in O(2^m) \approx O(1.3763^d)$.
\end{theorem}
\textbf{Proof}: We define $K$ and $M$ by $K = \frac{\log_2 1.5}{\log_2
  4.5} (d-1)$ and $M = \frac{1}{\log_2 4.5} (d-1)$.  We note that
\[
2K+M = \frac{2 \log_2 1.5 + 1}{\log_2 4.5} (d-1) = \frac{\log_2 \left(
1.5^2 \cdot 2 \right)}{\log_2 4.5} (d-1) = d-1.
\]
Since $k = \lceil K \rceil$, we have $K \leq k < K+1$, and since $2k+m
= 2K+M$, this implies $M-2 < m \leq M$.  Furthermore, $K = M \log_2
1.5$; therefore, $2^K = 1.5^M$, or equivalently, $\left( \frac{4}{3}
\right)^M 2^K = 2^M$.  In particular, $\Theta \left( \left(
\frac{4}{3} \right)^m 2^k \right) = \Theta \left( \left( \frac{4}{3}
\right)^M 2^K \right) = \Theta(2^M) = \Theta(2^m)$.  Thus, $\left(
\frac{4}{3} \right)^m (2^{k+1}) + 2^{m+1} \in \Theta\left( \left(
\frac{4}{3} \right)^m 2^k \right)$.  From Corollary~\ref{splitting
antipodal}, this implies that $|D_m| \in O(2^m) = O(2^M) \approx
O(1.3763^d)$.\Bx

\bibliographystyle{plain}

\end{document}